\documentclass[12pt,leqno]{article} 
\usepackage[latin1]{inputenc}
\usepackage[T1]{fontenc}
\usepackage[english]{babel}
\usepackage{times}
\usepackage{amssymb} 
\usepackage{graphicx,subfigure} 
\usepackage{tikz}
\usepackage{bm,dsfont}
\usepackage[margin=1.8cm,top=2.8cm]{geometry}
\usepackage{theorem}
\usepackage{multicol}
\usepackage{float}
\usepackage{hyperref}

\begin{document}
\title{
\textbf{On a nonlocal Boussinesq system \\ for internal wave propagation}}

\author{A.~Dur\'an\thanks{Departamento de Matem\'{a}tica Aplicada, Universidad Universidad de Valladolid,
Paseo Bel\'en 15, 47011--Valladolid(SPAIN).
Email: angel@mac.uva.es},}

\markboth{~\hrulefill\ A. Dur\'an}{On a nonlocal Boussinesq system for internal wave propagation\hrulefill~}

\maketitle

\abstract{
In this paper we are concerned with a nonlocal system to model the propagation of internal waves in a two-layer interface problem with rigid lid assumption and under a Boussinesq regime for both fluids. The main goal is to investigate aspects of well-posedness of the Cauchy problem for the deviation of the interface and the velocity, as well as the existence of solitary wave solutions and some of their properties.
}


\section{Introduction}
\label{sec:sec0}
\begin{figure}[htbp]
\centering
{\includegraphics[width=0.9\textwidth]{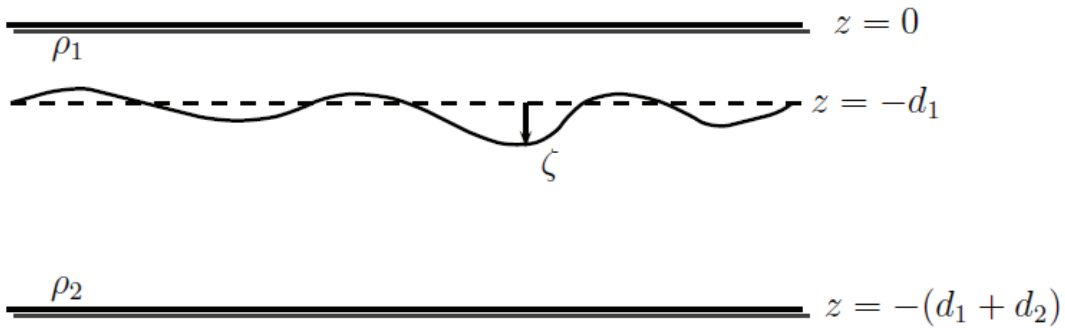}}
\caption{Idealized model of internal wave propagation in a two-layer interface.}
\label{f0}
\end{figure}

In this work a one-dimensional, nonlocal differential system for internal waves is considered. The system is derived in \cite{BLS2008} and describes the propagation of internal waves in a two-layer interface problem with rigid lid assumption and under the  Boussinesq regime for both fluids. The idealized model is sketched in Figure \ref{f0}. 

This consists of two inviscid, homogeneous, incompressible fluids of depths $d_{j}, j=1,2$ and densities $\rho_{j}, j=1,2$ with $\rho_{2}>\rho_{1}$. The system is bounded above and below by a rigid horizontal plane, with the origin of the vertical coordinate $z$ at the top
The deviation of the interface, denoted by $\zeta$, is assumed to be a graph over the bottom (described by the variable $x$) and surface tension effects are not considered.
The approach in \cite{BLS2008} is based on the reformulation of the Euler system with two nonlocal operators that link the velocity potentials associated to the layers at the interface. Then, by using suitable asymptotic expansions of these operators, several asymptotic models, consistent with the Euler system, are derived. They are associated to different physical regimes for the layers. The one considered here is the so-called Boussinesq-Boussinesq (B/B) regime, for which the interfacial deformations are assumed to be of small amplitude for both the upper and lower fluid domains and, additionally, the flow has a Boussinesq structure with respect to the two layers, with the nonlinear and dispersive effects of the same size for both fluids, \cite{BChS2002,BChS2004,BLS2008}.

One of the differential systems to model the B/B regime in 1D is given by
\begin{eqnarray}
\zeta_{t}+\frac{1}{\gamma+\delta}\partial_{x}{v}_{\beta}
+\left(\frac{\delta^{2}-\gamma}{(\delta+\gamma)^{2}}\right)\partial_{x} \left(\zeta{ v}_{\beta}\right)=0,&&\nonumber\\
(1-\beta\partial_{xx})({v}_{\beta})_{t}+(1-\gamma)\partial_{x}\zeta+\left(\frac{\delta^{2}-\gamma}{2(\delta+\gamma)^{2}}\right)\partial_{x} {v}_{\beta}^{2}=0,&&\label{bbs1}
\end{eqnarray}
where $\gamma=\rho_{1}/\rho_{2}<1, \delta=d_{1}/d_{2}$ are, respectively, the density and depth ratios, and $\beta=\displaystyle\frac{(1+\gamma\delta)}{3\delta(\gamma+\delta)}$. The variables $x$ and $t$ are proportional to distance along the channel and time respectively and $u(x,t)$ is the velocity variable with $v_{\beta}=(1-\beta\partial_{xx})^{-1}u$. When $\gamma=0, \delta=1$, (\ref{bbs1}) reduces to the classical Boussinesq system for surface water waves, see e.~g. \cite{Bou1871,BChS2002,BChS2004,PegoW1997,Whitham1974}. 
The purpose of this work is to deal with several properties of (\ref{bbs1}) as approximation of the Euler equations for internal waves:
\begin{enumerate}
\item The first point treated here is a review of the derivation of (\ref{bbs1}) without using the nonlocal operators defined and considered in \cite{BLS2008} but with the original variables instead, \cite{GrimshawP1998}.
\item Then, well-posedness of the corresponding Cauchy problem is considered. After the analysis of the linear case (cf. \cite{BLS2008}), the strategy used in \cite{BChS2002,BChS2004} for the surface wave problem is applied here to obtain a result of local existence and uniqueness of solution in suitable Sobolev spaces.
\item A final point of the paper is devoted to the solitary wave solutions of (\ref{bbs1}), that is, wave profiles traveling with permanent form and constant speed, decaying to zero at infinity. The study is divided into two parts. The first one is theoretical and proves the existence of such solutions for a range of speeds depending on the depth and density ratios of the two-fluid system, with the techniques considered in \cite{PegoW1997} for the case of surface waves. Since the theoretical study does not provide, in general, exact formulas for the waves, the second part is computational. A numerical technique, based on the Petviashvili iteration, \cite{Petv1976}, along with extrapolation, \cite{smithfs}, is applied to analyze, by numerical means, some properties of the solitary wave profiles, mainly focused on the regularity of the waves, their asymptotic decay and the speed-amplitude relation.
\end{enumerate}
All these results are displayed in the paper according to the following structure. In Section \ref{sec:sec2}, an alternative to \cite{BLS2008} for the derivation of (\ref{bbs1}) directly from the Euler equations, is provided. Section \ref{sec:sec3} is concerned with the analysis of well-posedness of the Cauchy problem. In Section \ref{sec:sec4} the existence result for solitary wave solutions is derived, while in Section \ref{sec:sec5} the computational study of the waves is carried out. We summarize the conclusions in Section \ref{sec:sec6}.


The following notation will be used throughout the paper. From the Sobolev space $H^{s}=H^{s}(\mathbb{R}), s\geq 0$ (with $H^{0}=L^{2}(\mathbb{R})$, the space of squared integrable functions on $\mathbb{R}$) we consider the product $X_{s_{1},s_{2}}=H^{s_{1}}\times H^{s_{2}}, s_{1},s_{2}\geq 0$, with associated norm
$$||u||_{s_{1},s_{2}}=\left(||u_{1}||_{s_{1}}^{2}+||u_{2}||_{s_{2}}^{2}\right)^{1/2},\quad u=(u_{1},u_{2})\in X_{s_{1},s_{2}},$$ where $||\cdot ||_{s}$ stands for the norm in $H^{s}$ and for simplicity we write $X_{s}:=X_{s,s}, s\geq 0$. For $T>0$ and $s\geq 0$, the space of continuous functions $v:(0,T]\rightarrow H^{s}$ is denoted by $C(0,T,H^{s})$ and its norm by
$$||v||_{C(0,T,H^{s})}=\max_{0<t\leq T}||v(t)||_{s}.$$ For $s_{1},s_{2}\geq 0$, the corresponding product space, $X_{T}^{s_{1},s_{2}}=C(0,T,H^{s_{1}})\times C(0,T,H^{s_{2}})$, is provided by the norm
$$||(v,w)||_{X_{T}^{s_{1},s_{2}}}=\left(||v||_{C(0,T,H^{s_{1}})}^{2}+||w||_{C(0,T,H^{s_{2}})}^{2}\right)^{1/2},$$ and where $X_{T}^{s}:=X_{T}^{s,s}, s\geq 0$. We will additionally make use of the Fourier transform
$$\widehat{f}(k)=\int_{\mathbb{R}}f(x)e^{-ikx}dx,\quad k\in\mathbb{R},\quad f\in H^{0}.$$
\section{Derivation of the model}
\label{sec:sec2}
In this section we shall derive the system (\ref{bbs1}) without using the approach based on nonlocal operators of \cite{BLS2008}. 
\subsection{Euler system for internal waves}
Assuming that each flow is irrotational, let $\Phi_{i}, i=1,2$ be the velocity potential associated to the upper and lower fluid layer respectively . For the model idealized in Figure \ref{f0}, the Euler equations can be written in dimensional, unscaled variables as follows. For $t>0$, the governing equations are
\begin{eqnarray}
\Phi_{1xx}+\Phi_{1zz}&=&0,\quad (x,z)\in\Omega_{t}^{1},\label{l21}\\
\Phi_{2xx}+\Phi_{2zz}&=&0,\quad (x,z)\in\Omega_{t}^{2},\label{l22}
\end{eqnarray}
where
\begin{eqnarray*}
\Omega_{t}^{1}&=&\{(x,z)/ -\infty<x<\infty, -d_{1}+\zeta(x,t)<z<0\},\\
\Omega_{t}^{2}&=&\{(x,z)/ -\infty<x<\infty, -d_{1}-d_{2}<z<-d_{1}+\zeta(x,t)\}.
\end{eqnarray*}
The rigid lid assumptions mean that the normal component of both velocity potentials vanishes at the corresponding boundary; that is, for $t>0$,
\begin{eqnarray}
\Phi_{1z}&=&0,\quad {\rm at}\quad \Gamma_{1}=\{(x,z)/ -\infty<x<\infty, z=0\}\label{l23}\\
\Phi_{2z}&=&0,\quad {\rm at}\quad \Gamma_{2}=\{(x,z)/ -\infty<x<\infty, z=-(d_{1}+d_{2})\}\label{l24}
\end{eqnarray}
The kinematic conditions at the interace between the fluids are
\begin{eqnarray}
\zeta_{t}+\Phi_{ix}\zeta_{x}=\Phi_{iz}\quad {\rm at}\quad \Gamma_{t}=\{(x,z)/ -\infty<x<\infty, z=-d_{1}+\zeta(x,t)\},\label{l25}
\end{eqnarray}
for $t>0$ and $i=1,2$. Finally, the condition of continuity of the pressure can be written as 
\begin{eqnarray}
\rho_{1}\left(\Phi_{1t}+\frac{1}{2}(\Phi_{1x}^{2}+\Phi_{1z}^{2})+g\zeta\right)=
\rho_{2}\left(\Phi_{2t}+\frac{1}{2}(\Phi_{2x}^{2}+\Phi_{2z}^{2})+g\zeta\right)\quad {\rm at} \quad \Gamma_{t},\label{l26}
\end{eqnarray}
where $g$ is the acceleration of gravity. From (\ref{l25}) and (\ref{l26}), two additional equations will be used throughout the section. The first one is obtained substracting the equations of (\ref{l25}) and eliminating $\zeta_{t}$:
\begin{eqnarray}
\Phi_{1x}\zeta_{x}-\Phi_{1z}=\Phi_{2x}\zeta_{x}-\Phi_{2z}
\quad {\rm at}\quad \Gamma_{t}=\{(x,z)/ -\infty<x<\infty, -z=-d_{1}+\zeta(x,t)\},\label{l25b}
\end{eqnarray}
while (\ref{l26}) can also be written as
\begin{eqnarray}
\Phi_{2t}-\gamma\Phi_{1t}+g(1-\gamma)\zeta +\frac{1}{2}(\Phi_{2x}^{2}-\gamma\Phi_{1x}^{2}+\Phi_{2z}^{2}-\gamma\Phi_{1z}^{2})
=0\quad {\rm at} \quad \Gamma_{t},\label{l26b}
\end{eqnarray}
\subsection{Non-dimensionalization}
In \cite{BLS2008}, the system (\ref{l21})-(\ref{l26}) is rewritten in terms of the potentials at the interface and two nonlocal operators which relate them. Instead, we directly nondimensionalize (\ref{l21})-(\ref{l26}) by using the dimensionless parameters
\begin{eqnarray*}
\epsilon=\frac{a}{d_{1}},\quad \mu=\left(\frac{d_{1}}{\lambda}\right)^{2},
\end{eqnarray*}
where $a$ is a typical amplitude and $\lambda$ a typical wavelength of the waves, as well as the dimensionless variables and unknowns
\begin{eqnarray*}
x=\lambda\widetilde{x},\quad z=d_{z}\widetilde{z},\quad t=\frac{\lambda}{\sqrt{gd_{1}}}\widetilde{t},\quad
\zeta=a\widetilde{\zeta},\quad \Phi_{i}=a\lambda\sqrt{\frac{g}{d_{1}}}\widetilde{\Phi_{i}}.
\end{eqnarray*}
Then the regions and boundaries are transformed into
\begin{eqnarray*}
\Omega_{t}^{1}&=&\{(x,z)/ -\infty<x<\infty, -1+\epsilon\zeta(x,t)<z<0\},\\
\Omega_{t}^{2}&=&\{(x,z)/ -\infty<x<\infty, -d_{1}-d_{2}<z<-d_{1}+\zeta(x,t)\},\\
\Gamma_{1}&=&\{(x,z)/ -\infty<x<\infty, z=0\},\quad 
 \Gamma_{2}=\{(x,z)/ -\infty<x<\infty, z= -1-\frac{1}{\delta}\},\\
 \Gamma_{t}&=&\{(x,z)/ -\infty<x<\infty, z=-1+\epsilon\zeta(x,t)\},
\end{eqnarray*}
while (\ref{l21})-(\ref{l26}) (see also (\ref{l26b}) are written as
\begin{eqnarray}
\mu\Phi_{1xx}+\Phi_{1zz}&=&0,\quad (x,z)\in\Omega_{t}^{1},\label{l21a}\\
\mu\Phi_{2xx}+\Phi_{2zz}&=&0,\quad (x,z)\in\Omega_{t}^{2},\label{l22a}\\
\Phi_{1z}&=&0,\quad {\rm at}\quad \Gamma_{1},\label{l23a}\\
\Phi_{2z}&=&0,\quad {\rm at}\quad \Gamma_{2},\label{l24a}\\
\zeta_{t}+\epsilon\Phi_{ix}\zeta_{x}&=&\frac{1}{\mu}\Phi_{iz}\quad {\rm at}\quad \Gamma_{t},\label{l25a}\\
\Phi_{2t}-\gamma\Phi_{1t}+(1-\gamma)\zeta +\frac{\epsilon}{2}(\Phi_{2x}^{2}-\gamma\Phi_{1x}^{2})+\frac{\epsilon}{2\mu}(\Phi_{2z}^{2}-\gamma\Phi_{1z}^{2})
&=&0\quad {\rm at}\quad  \Gamma_{t},\label{l26a}
\end{eqnarray}
where $t>0$. (Tildes were dropped.) 
\subsection{Boussinesq/Boussinesq regime}
We assume $\delta \sim 1$ and that the deviation of the interface is long and of small amplitude for both fluids, which means that $\epsilon,\mu<<1$ as well as $\epsilon_{2},\mu_{2}<<1$, where
\begin{eqnarray*}
\epsilon_{2}=\frac{a}{d_{2}}=\epsilon\delta,\quad \mu_{2}=\left(\frac{d_{2}}{\lambda}\right)^{2}=\frac{\mu}{\delta^{2}},
\end{eqnarray*}
are the amplitude and wavelength parameters with respect to the lower fluid (and which are not independent of $\epsilon, \mu$). Furthermore, we are interested in the so-called Boussines/Boussinesq regime, \cite{BLS2008}. This means that the nonlinear and dispersive effects ar of the same size for both fluids and thus
\begin{eqnarray}
\epsilon\sim\mu\sim\epsilon_{2}\sim\mu_{2}.\label{l27}
\end{eqnarray}
From this point, the derivation follows the strategy considered in \cite{BChS2002} for surface waves, see also \cite{GrimshawP1998,Peregrine1972,Benjamin1972,Whitham1974}. The potentials $\Phi_{i}, i=1,2$, are formally expanded in the corresponding domains as
\begin{eqnarray}
\Phi_{1}(x,z,t)&=&\sum_{m=0}^{\infty}f_{1m}(x,t)z^{m},\quad (x,z)\in\Omega_{t}^{1},\label{l281}\\
\Phi_{2}(x,z,t)&=&\sum_{m=0}^{\infty}f_{2m}(x,t)(z+h_{\delta})^{m},\quad (x,z)\in\Omega_{t}^{2},\label{l282}
\end{eqnarray}
where $h_{\delta}=1+\displaystyle\frac{1}{\delta}$. Now, satisfaction of (\ref{l21a}) and (\ref{l22a}) implies
\begin{eqnarray}
(m+2)(m+1)f_{i,m+2}(x,t)=-\mu(f_{im}(x,t))_{xx},\quad i=1,2,\label{l28}
\end{eqnarray}
while, due to the boundary conditions (\ref{l23a}) and (\ref{l24a}), we have $f_{i,1}(x,t)=0, i=1,2$, and therefore, according to (\ref{l28}), $f_{i,2k+1}(x,t)=0, k\geq 0, i=1,2$. For the even terms, the application of (\ref{l28}) leads to
\begin{eqnarray*}
f_{i,2k}(x,t)=\frac{(-1)^{k}}{(2k)!}\mu^{k}\frac{\partial^{2k}}{\partial x^{2k}}F_{i}(x,t),\quad k\geq 0, i=1,2,
\end{eqnarray*}
where $F_{i}(x,t)=f_{i,0}(x,t), i=1,2$. Therefore, (\ref{l281}), (\ref{l282}) can be written as
\begin{eqnarray}
\Phi_{1}(x,z,t)&=&\sum_{k=0}^{\infty}\frac{(-1)^{k}}{(2k)!}\mu^{k}\frac{\partial^{2k}}{\partial x^{2k}}F_{1}(x,t)z^{m},\quad (x,z)\in\Omega_{t}^{1},\label{l291}\\
\Phi_{2}(x,z,t)&=&\sum_{m=0}^{\infty}\frac{(-1)^{k}}{(2k)!}\mu^{k}\frac{\partial^{2k}}{\partial x^{2k}}F_{2}(x,t)(z+h_{\delta})^{m},\quad (x,z)\in\Omega_{t}^{2},\label{l292}
\end{eqnarray}
In what follows, the expansions (\ref{l291}) and (\ref{l292}) will be evaluated at the interface $z=-1+\epsilon\zeta$ and, according to (\ref{l27}), the linear terms in $\epsilon, \mu$ will be retained. Note first that (\ref{l25b}), in nondimensional form, reads
\begin{eqnarray}
\epsilon\Phi_{1x}\zeta_{x}-\frac{1}{\mu}\Phi_{1z}&=&\epsilon\Phi_{2x}\zeta_{x}-\frac{1}{\mu}\Phi_{2z}\label{l25bb}\\
&& {\rm at}\quad \Gamma_{t}=\{(x,z)/ -\infty<x<\infty, z=-1+\epsilon\zeta(x,t)\},\nonumber
\end{eqnarray}
which, in terms of the expansions (\ref{l291}) and (\ref{l292}), leads to
\begin{eqnarray}
\frac{1}{\delta}F_{2x}+F_{1x}=\epsilon\zeta(F_{1x}-F_{2x})+\frac{\mu}{6}(F_{1xxx}+\frac{1}{\delta^{3}}F_{2xxx})+O(\epsilon^{2},\epsilon\mu,\mu^{2}).\label{l210}
\end{eqnarray}
On the other hand, we define $\Phi=\Phi_{2}-\gamma\Phi_{1}$ and let $u=\Phi_{x}$. According to (\ref{l291}) and (\ref{l292}), we have 
\begin{eqnarray}
u=F_{2x}-\gamma F_{1x}-\frac{\mu}{2}\left(\frac{1}{\delta^{2}}F_{2xxx}-\gamma F_{1xxx}\right)+O(\epsilon^{2},\epsilon\mu,\mu^{2}).\label{l211}
\end{eqnarray}
Formulas (\ref{l210}) and (\ref{l211}) are now used to determine iteratively $F_{1}$ and $F_{2}$ in terms of $\zeta$ and $u$ (cf. \cite{GrimshawP1998}). After some computations, we have
\begin{eqnarray}
F_{1x}&=&\frac{1}{\delta+\gamma}\left(-u-\epsilon\frac{\delta(1+\delta)}{\delta+\gamma}\zeta u-\mu\frac{2+3\gamma\delta+\delta^{2}}{6\delta (\delta+\gamma)}u_{xx}\right)+O(\epsilon^{2},\epsilon\mu,\mu^{2}).\label{l212a}\\
F_{2x}&=&\frac{1}{\delta+\gamma}\left(\delta u-\epsilon\frac{\delta\gamma(1+\delta)}{\delta+\gamma}\zeta u+\mu\frac{1+3\delta+\delta^{2}(3\gamma-1)}{6\delta (\delta+\gamma)}u_{xx}\right)+O(\epsilon^{2},\epsilon\mu,\mu^{2}).\label{l212b}
\end{eqnarray}
Finally, the application of (\ref{l212a}) and (\ref{l212b}) to the first kinematic condition in (\ref{l25a}) and, after differentiating in $x$, to (\ref{l26a}), imply
\begin{eqnarray}
\zeta_{t}+\frac{1}{\delta+\gamma}\partial_{x}(1+\beta\mu\partial_{xx})u+\epsilon\left(\frac{\delta^{2}-\gamma}{(\delta+\gamma)^{2}}\right)\partial_{x}(\zeta u)&=&O(\epsilon^{2},\epsilon\mu,\mu^{2}),\label{l213a}\\
u_{t}+(1-\gamma)\zeta_{x}+\frac{\epsilon}{2}\left(\frac{\delta^{2}-\gamma}{(\delta+\gamma)^{2}}\right)\partial_{x}(u^{2})&=&O(\epsilon^{2},\epsilon\mu,\mu^{2}),\label{l213b}
\end{eqnarray}
where  $\beta$ is defined in (\ref{bbs1}). Note that if 
$v_{\beta}=(1-\beta\mu\partial_{xx})^{-1}u$, then $v_{\beta}=(1+\beta\mu\partial_{xx})u+O(\mu^{2})$ and substitution into  (\ref{l213a}) and (\ref{l213b}) leads to (\ref{bbs1}), assuming (\ref{l27}) and when the $O(\epsilon^{2},\epsilon\mu,\mu^{2})$ terms are dropped.

{\bf Remark.}
From system (\ref{l213a}), (\ref{l213b}), the three parameter family of Boussinesq systems for internal waves presented in \cite{BLS2008} can alternatively be derived by adapting the procedure implemented in \cite{BChS2002} for the case of surface waves.

\section{Well-posedness} \label{sec:sec3}

\subsection{Linear well-posedness}
We now study the linear well-posedness of (\ref{bbs1}). The associated linear system, written in terms of $\zeta$ and $u$ is of the form
\begin{eqnarray}
&&\zeta_{t}+\frac{1}{\gamma+\delta}\partial_{x} (I-\beta\partial_{xx})^{-1}{u}_{x}=0,\nonumber\\
&&{u}_{t}+(1-\gamma)\partial_{x}\zeta=0.\label{bbs2}
\end{eqnarray}
The application of the Fourier transform leads to the system
\begin{eqnarray*}
\frac{d}{dt}\left(\matrix{\widehat{\zeta}(k,t)\cr \widehat{u}(k,t)}\right)+(ik)A(k)\left(\matrix{\widehat{\zeta}(k,t)\cr \widehat{u}(k,t)}\right)=0,
\quad k\in\mathbb{R},
\end{eqnarray*}
where
\begin{eqnarray}
A(k)=\left(\matrix{0&\omega(k)\cr 1-\gamma&0}\right),\quad 
\omega(k)=\frac{1}{\delta+\gamma}\frac{1}{1+\beta k^{2}},\label{bbs4}
\end{eqnarray}
Then well-posedness is determined by the behaviour of the matrix
\begin{eqnarray}
e^{-ikA(k)t}=\left(\matrix{\cos(k\sigma(k)t)&-i\sqrt{\frac{\omega(k)}{1-\gamma}}\sin(k\sigma(k)t)\cr i\sqrt{\frac{1-\gamma}{\omega(k)}}\sin(k\sigma(k)t)& \cos(k\sigma(k)t)}\right),
\label{bbs3b}
\end{eqnarray}
where $\sigma(k)=\sqrt{(1-\gamma)\omega(k)}$. Specifically, the linearized problem (\ref{bbs2}) is well-posed when (\ref{bbs3b}) is bounded in finite intervals of $t$. We define the order of $\sigma(k)$ as the integer $l$ such that
$$\sigma(k)\approx |k|^{l}, \quad |k|\rightarrow\infty,$$  and let $m_{1}=\max\{0,-l\}, m_{2}=\max\{0,l\}$. Since $\omega(k)$ has neither poles or zeros on the real axis then we have, \cite{BChS2002}:

{\bf Theorem 3.1.}
{\em 
Let $\beta>0$. System (\ref{bbs1}) is linearly well-posed for $(\zeta,u)\in X_{s+1,s}$ and therefore for $(\zeta,v_{\beta})\in X_{s+1,s+2}, s\geq 0$.
}

{\bf Remark.} Theorem 3.1 completes the linear well-posedness result established in \cite{BLS2008} by specifying the spaces where that holds.
\subsection{Local well-posedness of the full nonlinear system}
Taking the Fourier transform in $x$ to (\ref{bbs1}) for $\zeta$ and $u$ we have 
\begin{eqnarray}
\frac{d}{dt}\left(\matrix{\widehat{\zeta}\cr\widehat{u}}\right)+i{k}{A}({ k})\left(\matrix{\widehat{\zeta}\cr\widehat{u}}\right)+ik{ F}\left(\matrix{\widehat{\zeta}\cr\widehat{u}}\right)(k)=0,\label{bbs5}
\end{eqnarray}
where $A$ is given by (\ref{bbs4}) and
\begin{eqnarray*}
{F}=K_{\gamma,\delta}\left(\matrix{\left[\zeta (1-\beta\partial_{xx})^{-1}u\right]^{\widehat{}}(k)\cr\frac{\left(\left[ (1-\beta\partial_{xx})^{-1}u)\right]^{2}\right)^{\widehat{}}}{2}(k)}\right),\quad 
K_{\gamma,\delta}=\displaystyle\frac{\delta^{2}-\gamma}{(\delta+\gamma)^{2}}.
\end{eqnarray*}
In order to study well-posedness of the system (\ref{bbs1}), we use a
similar strategy to that of \cite{BChS2004} by decoupling the linear part with the operator $\Sigma$ such that
\begin{eqnarray*}
\widehat{\Sigma f}(k)=\sqrt{\frac{\omega(k)}{1-\gamma}}\widehat{f}(k), \quad k\in\mathbb{R},
\end{eqnarray*}
where $\omega$ is given in (\ref{bbs4}). Defining new variables $v, w$ such that
\begin{eqnarray*}
\zeta=\Sigma(v+w),\quad u=v-w,
\end{eqnarray*}
then the system (\ref{bbs5}) is of the form
\begin{eqnarray}
\frac{d}{dt}\left(\matrix{\widehat{v}\cr\widehat{w}}\right)+i{k}\left(\matrix{\sigma(k)&0\cr 0&-\sigma(k)}\right)\left(\matrix{\widehat{v}\cr\widehat{w}}\right)+ikP^{-1}{ F}=0,&&\nonumber\\
&&\label{bbs6}
\end{eqnarray}
with $\sigma(k)$ given in (\ref{bbs3b}) and
\begin{eqnarray}
P^{-1}=\frac{1}{2}\left(\matrix{\sqrt{\frac{1-\gamma}{\omega(k)}}&1\cr \sqrt{\frac{1-\gamma}{\omega(k)}}&-1}\right).\label{bbs7}
\end{eqnarray}
Then the following local well-posedness result holds:

{\bf Theorem 3.2.}
{\em 
Let $(\zeta_{0},v_{\beta,0})\in X_{s+1,s+2}, s\geq 0$. Then there exists $T>0$ and a unique solution $(\zeta,u)$ in $X_{T}^{s+1,s}$ of (\ref{bbs1}) with initial condition $(\zeta_{0},u_{0})$ with $u_{0}=(1-\beta\partial_{xx})v_{\beta,0}$.
}

{\em Proof}. We consider the variables
\begin{eqnarray*}
v_{0}=\frac{\Sigma^{-1}(\zeta_{0})+u_{0}}{2},\quad 
w_{0}=\frac{\Sigma^{-1}(\zeta_{0})-u_{0}}{2}.
\end{eqnarray*}
Then $(v_{0},w_{0})\in X_{s}$ and taking the inverse Fourier transform in (\ref{bbs6}) we have
\begin{eqnarray}
\frac{d}{dt}\left(\matrix{{v}\cr{w}}\right)+
\mathcal{B}\left(\matrix{{v}\cr{w}}\right)=\mathcal{F}\left(\matrix{{v}\cr{w}}\right),\label{bbs8}
\end{eqnarray}
where $\mathcal{B}$ is the operator with symbol
\begin{eqnarray*}
i{k}\left(\matrix{\sigma(k)&0\cr 0&-\sigma(k)}\right),
\end{eqnarray*}
and 
\begin{eqnarray*}
\mathcal{F}\left(\matrix{{v}\cr{w}}\right)=-\mathcal{P}^{-1}K_{\gamma,\delta}\left(\matrix{[(1-\beta\partial_{xx})^{-1}(v-w)]\Sigma(v+w)\cr\frac{\left(\left[ (1-\beta\partial_{xx})^{-1}(v-w))\right]^{2}\right)}{2}}\right),
\end{eqnarray*}
with $\mathcal{P}$ the operator associated to $ik P(k)$, obtained from (\ref{bbs7}), as symbol. By Duhamel formula, the solution of (\ref{bbs8}) satisfies
\begin{eqnarray*}
\left(\matrix{{v}\cr{w}}\right)=\mathcal{S}(t)\left(\matrix{{v_{0}}\cr{w_{0}}}\right)+
\int_{0}^{t}\mathcal{S}(t-\tau)\mathcal{F}\left(\matrix{{v}\cr{w}}\right)d\tau,
\end{eqnarray*}
where $\mathcal{S}(t)$ is the group generated by $\mathcal{B}$.
Consider then the mapping $(\widetilde{v},\widetilde{w})\mapsto (v,w)$ where
\begin{eqnarray}
\qquad \left(\matrix{{v}\cr{w}}\right)=\mathcal{S}(t)\left(\matrix{{v_{0}}\cr{w_{0}}}\right)+
\int_{0}^{t}\mathcal{S}(t-\tau)\mathcal{F}\left(\matrix{\widetilde{v}\cr\widetilde{w}}\right)d\tau.\label{bbs9}
\end{eqnarray}
Note first that $\mathcal{S}(t)$ is a unitary group on $X_{s}$. On the other hand, the operator $\mathcal{F}$ can be estimated by using Lemma 2.2 of \cite{BChS2004}, in such a way that if $T>0$ and 
 if $(\widetilde{v}_{1},\widetilde{w}_{1}), (\widetilde{v}_{2},\widetilde{w}_{2})$ are in a closed ball of radius $R$ centered at $0$ in $X_{T}^{s}$ then there is $C(R)>0$ for which
\begin{eqnarray}
\left|\left|\mathcal{F}(\widetilde{v}_{1},\widetilde{w}_{1})(\tau)-
\mathcal{F}(\widetilde{v}_{2},\widetilde{w}_{2})(\tau)\right|\right|_{X_{s}}\leq C(R)\left|\left|\left(\matrix{\widetilde{v}_{1}\cr \widetilde{w}_{1}}\right)(\tau)-\left(\matrix{\widetilde{v}_{2}\cr \widetilde{w}_{2}}\right)(\tau)\right|\right|_{X_{s}},\label{bbs10}
\end{eqnarray}
for any $0\leq \tau\leq T$. Now, (\ref{bbs10}) and the property $\mathcal{F}(0,0)=(0,0)$ imply that there is $T>0$, sufficiently small, such that (\ref{bbs9}) is a contraction of the closed ball into itself. Then the result follows by using the Contraction Mapping Theorem.$\Box$
\section{Solitary wave solutions} \label{sec:sec4}
In this section, the existence of solitary wave solutions of (\ref{bbs1}) is studied. They are solutions of traveling-wave form $\zeta=\zeta(x-c_{s}t), u=u(x-c_{s}t)$ (or $v_{\beta}=v_{\beta}(x-c_{s}t)$) for some speed $c_{s}\neq 0$ and where the profiles $\zeta=\zeta(X), u=u(X), X=x-c_{s}t$ (or $v_{\beta}=v_{\beta}(X)$) with $\zeta,u\rightarrow 0$ as $|X|\rightarrow \infty$  must satisfy
\begin{eqnarray}
\qquad
\left(\matrix{c_{s}&\frac{-1}{\delta+\gamma}\cr \gamma-1&c_{s}(1-\beta\partial_{xx})}\right)\left(\matrix{\zeta\cr v_{\beta}}\right)=K_{\gamma,\delta}\left(\matrix{\zeta v_{\beta}\cr \frac{v_{\beta}^{2}}{2}}\right)
\label{bbs11}.
\end{eqnarray}
Solving the first equation of (\ref{bbs11}) for $\zeta$ and substituting into the second one lead to
\begin{eqnarray}
v_{\beta}^{\prime\prime}-\frac{1}{\beta}v_{\beta}+G^{\prime}(v_{\beta})=0,\label{bbs12}
\end{eqnarray}
where $G(v)=\displaystyle\int_{0}^{v}\frac{1}{\beta c_{s}}g(z)dz$,  $c_{\gamma,\delta}=\displaystyle\sqrt{\frac{1-\gamma}{\delta+\gamma}}$ and
\begin{eqnarray*}
g(v)=\frac{\delta^{2}-\gamma}{2(\delta+\gamma)^{2}}v^{2}+\frac{c_{\gamma,\delta}^{2}v}{c_{s}-\frac{\delta^{2}-\gamma}{(\delta+\gamma)^{2}}v}.\label{bbs13}
\end{eqnarray*}
\subsection{Existence of solitary waves}
The analysis of (\ref{bbs12}) leads to the following result (cf. \cite{PegoW1997}).

{\bf Theorem 4.1.}
{\em Let $0<\gamma<1, \delta^{2}-\gamma\neq 0$ and $c_{s}^{2}-c_{\gamma,\delta}^{2}>0$. Then (\ref{bbs11}) has a unique solution $(\zeta_{s}(X),v_{\beta,s}(X))$ which is even and goes to zero as $|X|\rightarrow \infty$. The profiles $\zeta, v_{\beta}$ are of elevation when $\delta^{2}-\gamma>0$ and of depression when $\delta^{2}-\gamma<0$.
}

\begin{figure}[htbp]
\centering
\subfigure[]
{\includegraphics[height=8.35cm,width=8.35cm]{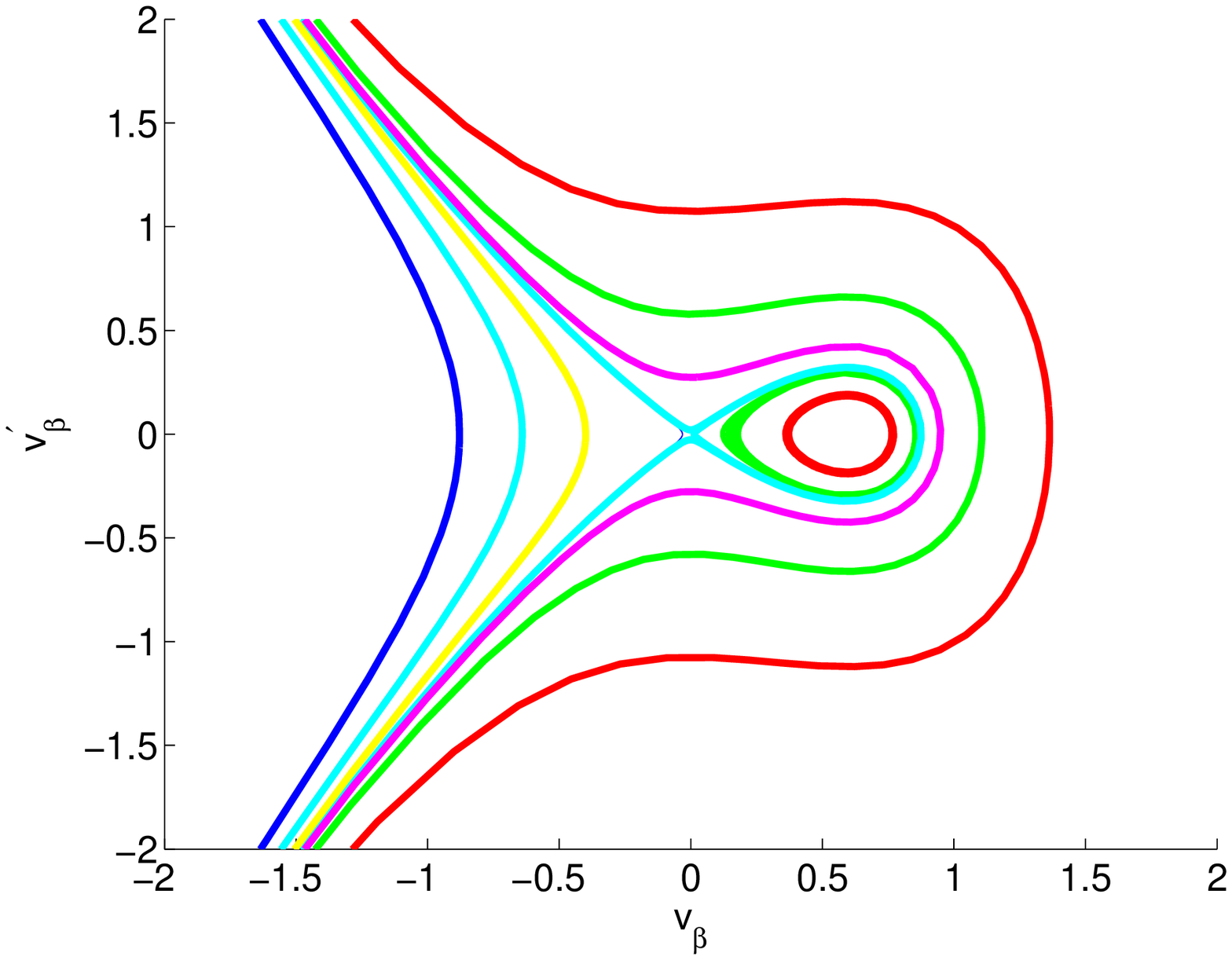}}
\subfigure[]
{\includegraphics[height=8.35cm,width=8.35cm]{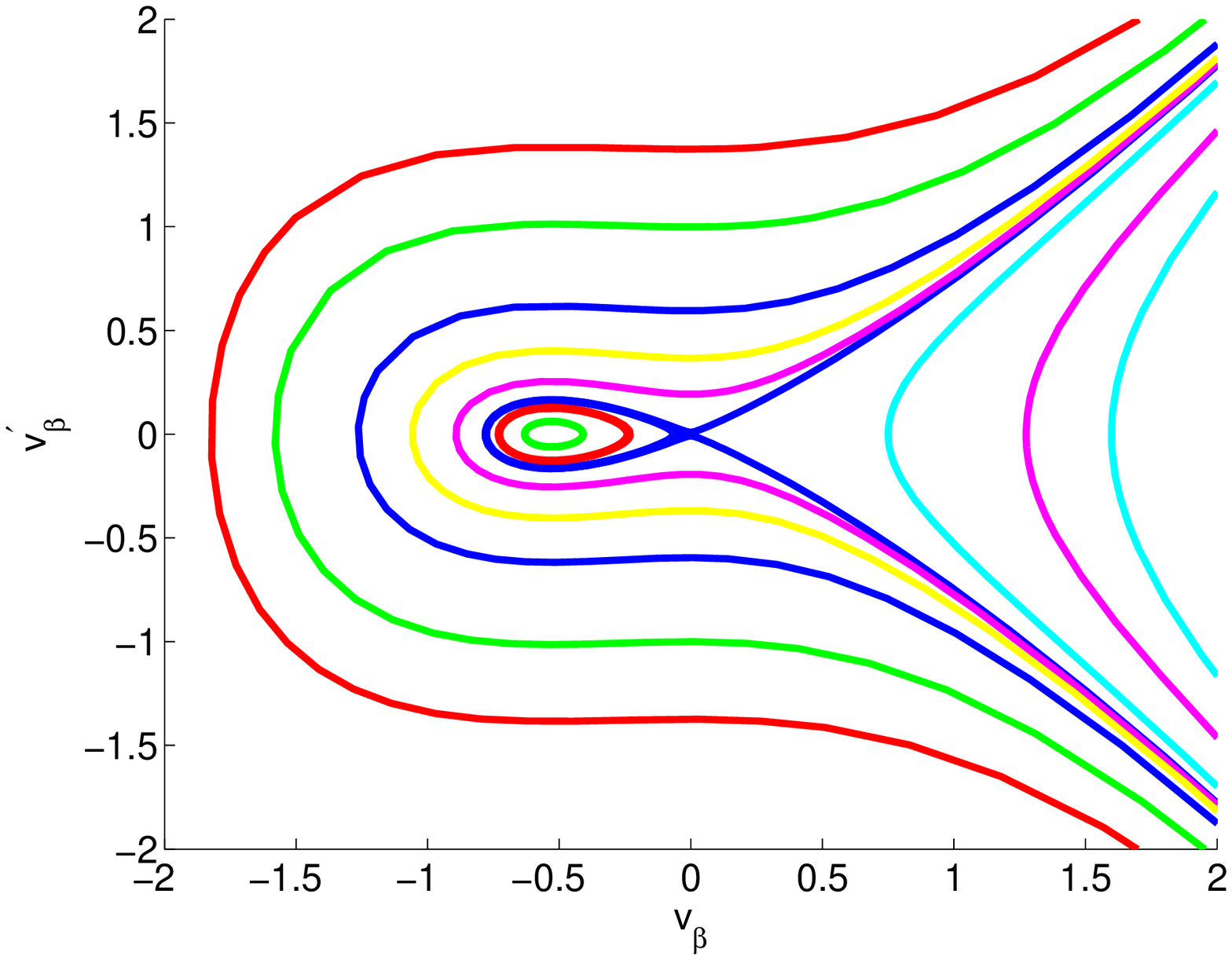}}
\caption{Phase portraits for (\ref{bbs12}) with (a) $K_{\gamma,\delta}>0$ and (b) $K_{\gamma,\delta}<0$.}
\label{f1}
\end{figure}

{\em Proof}. We assume $c_{s}>0$ (the arguments are similar in the case $c_{s}<0$). The equation (\ref{bbs12}) is conservative with the energy given by
$$E=\frac{1}{2}(v_{\beta}^{\prime})^{2}+U(v_{\beta}),$$ where 
$U(x)=-\displaystyle\frac{x^{2}}{2\beta}+G(x)$ is the potential energy. Since $U(0)=U^{\prime}(0)=0$ and $U^{\prime\prime}(0)=-\displaystyle\frac{c_{s}^{2}-c_{\gamma,\delta}^{2}}{\beta c_{s}^{2}}<0$ then the phase plane analysis shows that the origin is a saddle point (see Figure \ref{f1}). 
Note that
\begin{eqnarray*}
G(v)=\frac{K_{\gamma,\delta}v^{3}}{6\beta c_{s}}+\left(\frac{c_{s}}{c_{\gamma,\delta}}\right)^{2}\frac{c_{s}}{\beta K_{\gamma,\delta}}\left(-v+\frac{c_{s}}{K_{\gamma,\delta}}{\rm ln}\left(\frac{c_{s}/K_{\gamma,\delta}}{\frac{c_{s}}{K_{\gamma,\delta}}-v}\right)\right).
 \end{eqnarray*} 
 Note then that the sign of $K_{\gamma,\delta}$ determines the behaviour of $G$ in the following sense (cf. \cite{PegoW1997}):
 \begin{itemize}
 \item[(i)] If $K_{\gamma,\delta}>0$, then $G$ increases smoothly from $0$ to $\infty$ as $0<v<\displaystyle\frac{c_{s}}{K_{\gamma,\delta}}$.
 \item[(ii)]  If $K_{\gamma,\delta}<0$, then $G$ decreases smoothly from $\infty$ to $0$ as $\displaystyle\frac{c_{s}}{K_{\gamma,\delta}}<v<0$.
 \end{itemize}
 In both cases, the function $G$ is always positive and therefore we can find $v^{*}=v^{*}(\gamma,\delta)$ such that
 $$U(v^{*})=G(v^{*})-\frac{(v^{*})^{2}}{2\beta}=0,$$ with $0<v^{*}<\displaystyle\frac{c_{s}}{K_{\gamma,\delta}}$ in the case (i) and $\displaystyle\frac{c_{s}}{K_{\gamma,\delta}}<v^{*}<0$ in the case (ii). Thus, the solution of (\ref{bbs12}) with $v_{\beta}(0)=v^{*}, v_{\beta}^{\prime}(0)=0$ has zero energy, is even, positive when $K_{\gamma,\delta}>0$, negative when $K_{\gamma,\delta}<0$ and goes to zero as $|X|\rightarrow\infty$. Note that from (\ref{bbs11}) we have
 $$\zeta_{s}=\frac{1}{K_{\gamma,\delta}(\gamma+\delta)}\frac{v_{\beta}}{\frac{c_{s}}{K_{\gamma,\delta}}-v_{\beta}},$$ and from (\ref{bbs12})
 $$u_{s}:=(1-\beta \partial_{xx})v_{\beta}=\beta G^{\prime}(v_{\beta}).$$ Therefore, both $\zeta_{s}, u_{s}$ are also of elevation when $K_{\gamma,\delta}>0$ and of depression when $K_{\gamma,\delta}<0$.$\Box$
 \section{Numerical approximation to solitary waves}
 \label{sec:sec5}
 Since explicit formulas for the solitary waves are, to our knowledge, in general unknown, then the study of the form and additional properties of the profiles will be done in this section by computational means. A numerical procedure used to this end will be first briefly described and then applied to analyze features of the waves, mainly concerned with their regularity, asymptotic decay and the relation between the amplitude and the speed of the profiles.
 \subsection{A numerical technique to compute solitary wave profiles}
 In order to compute approximate solitary wave profiles of (\ref{bbs1}), the system (\ref{bbs11}) is discretized on an interval $(-l,l), l$ large, using Fourier pseudospectral approximation. Let $N\geq 1$ be integer and even, and let $(\zeta_{h},v_{h})$ be a $2N-$vector approximation to the profile $(\zeta,v_{\beta})$ at the uniform grid of collocation points $x_{j}=-l+jh, h=2l/N, j=0,\ldots,N$. Then $(\zeta_{h},v_{h})$ satisfies the system
 \begin{eqnarray}
&& \mathcal{L}_{h}\left(\matrix{\zeta_{h}\cr v_{h}}\right)=\mathcal{N}_{h}(\zeta_{h},v_{h})\label{l51}\\
&&\mathcal{L}_{h}=\left(\matrix{c_{s}I_{N}&-\frac{1}{\delta+\gamma}I_{N}\cr
(\gamma-1)I_{N}&c_{s}(1-\beta D_{N}^{2})}\right),\quad
\mathcal{N}_{h}(\zeta_{h},v_{h})=K_{\gamma,\delta}\left(\matrix{\zeta_{h}. v_{h}\cr \frac{1}{2}v_{h}.v_{h}}\right),\nonumber
 \end{eqnarray}
 where $I_{N}$ stands for the $N\times N$ identity matrix, $D_{N}$ stands for the $N\times N$ pseudospectral differentiation matrix, \cite{Canutohqz}, and the dot in $\mathcal{N}_{h}$ stands for the Hadamard product of vectors. System (\ref{l51}) is implemented in Fourier space, in such a way that if $\widehat{\zeta}_{h}(k),\widehat{v}_{h}(k), k=0,\ldots,N-1$ denote the $k-$th discrete Fourier components of $\zeta_{h}$ and $v_{h}$ respectively, and $k^{\prime}=\pi k/l$, then
 (\ref{l51}) leads to $2\times 2$ algebraic systems of the form
 \begin{eqnarray}
 \left(\matrix{c_{s}&-\frac{1}{\delta+\gamma}\cr
(\gamma-1)&c_{s}(1+\beta (k^{\prime})^{2})}\right)\left(\matrix{\widehat{\zeta}_{h}(k)\cr \widehat{v}_{h}(k)}\right)=K_{\gamma,\delta}\left(\matrix{\widehat{\zeta_{h}. v_{h}}(k)\cr \frac{1}{2}\widehat{v_{h}. v_{h}}(k)}\right),\quad k=0,\ldots,N-1.\label{l52}
 \end{eqnarray}
 Once (\ref{l52}) is solved for each $k=0,\ldots,N-1$the approximation $u_{h}$ to $u$ at the collocation points  is defined as $u_{h}=(I_{N}-\beta D_{N}^{2})v_{h}$ or, in Fourier components, as $\widehat{u_{h}}(k)=(1+\beta (k^{\prime})^{2})\widehat{v_{h}}(k), k=0,\ldots,N-1$.
 
 The system (\ref{l51}), or its Fourier version (\ref{l52}), is numerically solved by iteration with the Petviashvili method, \cite{Petv1976}: given initial data $\zeta_{h}^{[0]}, v_{h}^{[0]}$, the $(\nu+1)-$th iteration $\zeta_{h}^{[\nu+1]}, v_{h}^{[\nu+1]}, \nu=0,1,\ldots,$ is the solution of
 \begin{eqnarray}
 m_{\nu}&=&\frac{\langle  \mathcal{L}_{h}\left(\matrix{\zeta_{h}^{[\nu]}\cr v_{h}^{[\nu]}}\right),
 \left(\matrix{\zeta_{h}^{[\nu]}\cr v_{h}^{[\nu]}}\right)\rangle}{\langle  \mathcal{N}_{h}(\zeta_{h}^{[\nu]},v_{h}^{[\nu]}),
 \left(\matrix{\zeta_{h}^{[\nu]}\cr v_{h}^{[\nu]}}\right)\rangle},\nonumber\\
 \mathcal{L}_{h}\left(\matrix{\zeta_{h}^{[\nu+1]}\cr v_{h}^{[\nu+1]}}\right)&=&m_{\nu}^{2}\mathcal{N}_{h}(\zeta_{h}^{[\nu]},v_{h}^{[\nu]}),\label{l53}
 \end{eqnarray}
 where $\langle\cdot,\cdot\rangle$ stands for the Euclidean inner product in $\mathbb{R}^{2N}$. (See e.~g. \cite{PelinovskyS2004} for a justification of the formulas.) The iteration is supplemented with an extrapolation method, \cite{smithfs,jbilous}, a technique which has been revealed useful in the acceleration of the convergence when computing approximate solitary wave profiles, \cite{AlvarezD2015}. The implementation of the iterative procedure is carried out in Fourier space and the accuracy of the method was checked in the standard way, see e.~g. \cite{DDM1} for the details.
 \begin{figure}[htbp]
\centering
\subfigure[]
{\includegraphics[height=8.35cm,width=8.35cm]{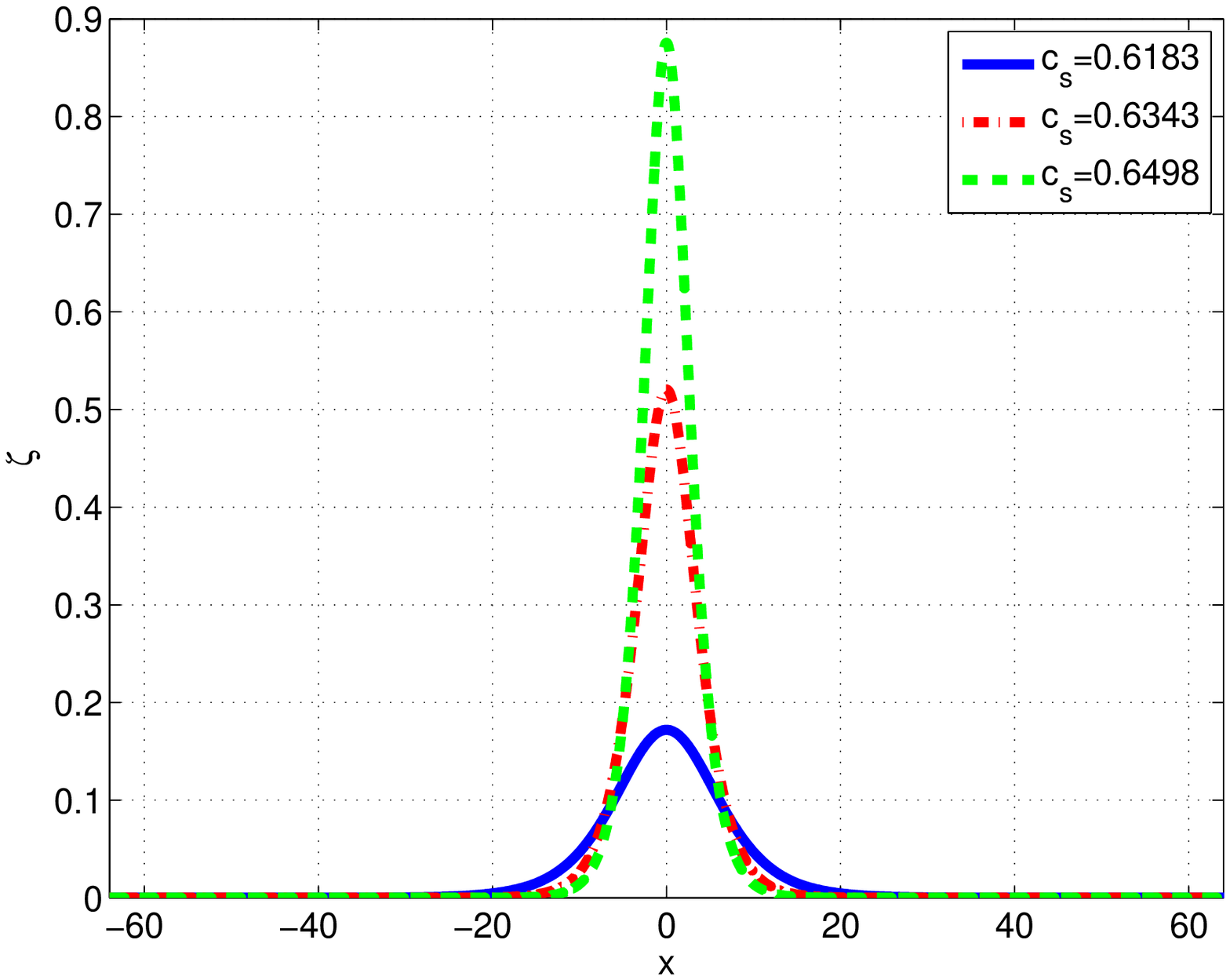}}
\subfigure[]
{\includegraphics[height=8.35cm,width=8.35cm]{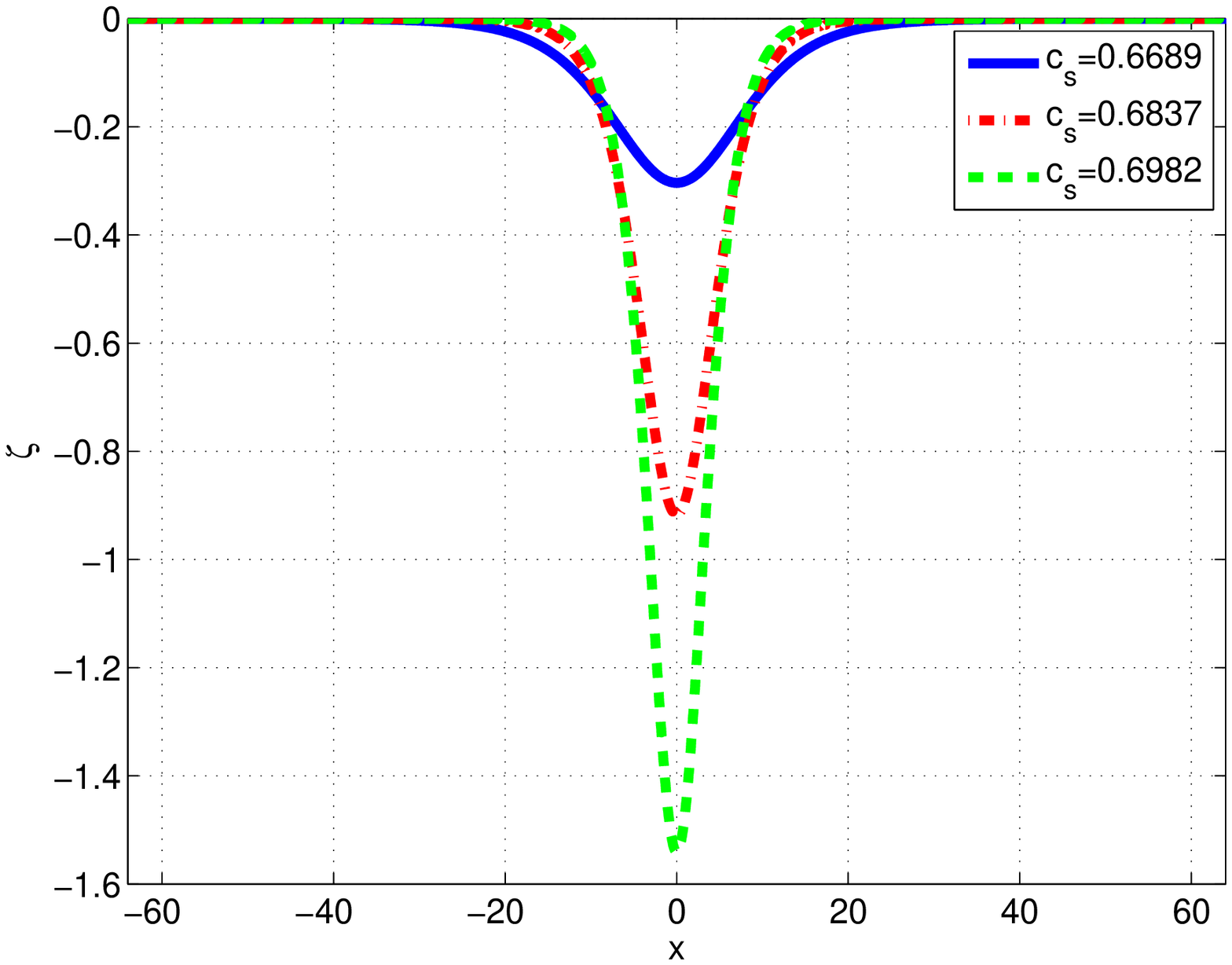}}
\caption{Approximate $\zeta$ solitary-wave profiles of (\ref{bbs11}): (a) corresponds to $\gamma=0.5, \delta=0.8$ ($K_{\gamma,\delta}>0$, waves of elevation) while (b)  corresponds to $\gamma= \delta=0.5$ ($K_{\gamma,\delta}<0$, waves of depression).}
\label{f2}
\end{figure}
\subsection{Numerical results} 
The form of the profiles of depression and of elevation is illustrated in Figure \ref{f2}, which displays some computed waves for different speeds, close to the limiting value $c_{\gamma,\delta}$. Note that the amplitude (the magnitude of the maximum negative excursion for the case of waves of depression) increases with the speed. 
\begin{figure}[htbp]
\centering
\subfigure[]
{\includegraphics[height=5.9cm,width=5.9cm]{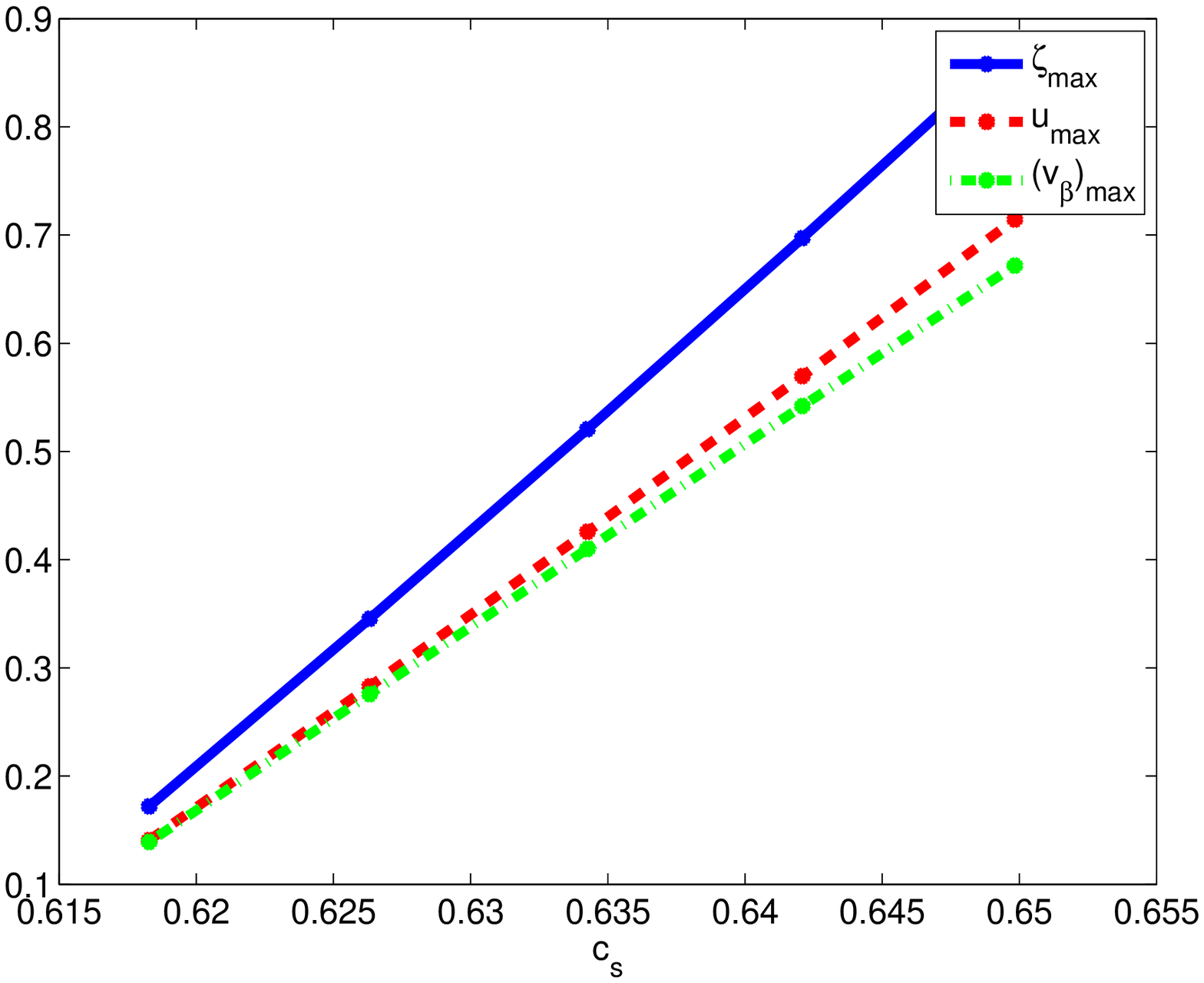}}
\subfigure[]
{\includegraphics[height=5.9cm,width=5.9cm]{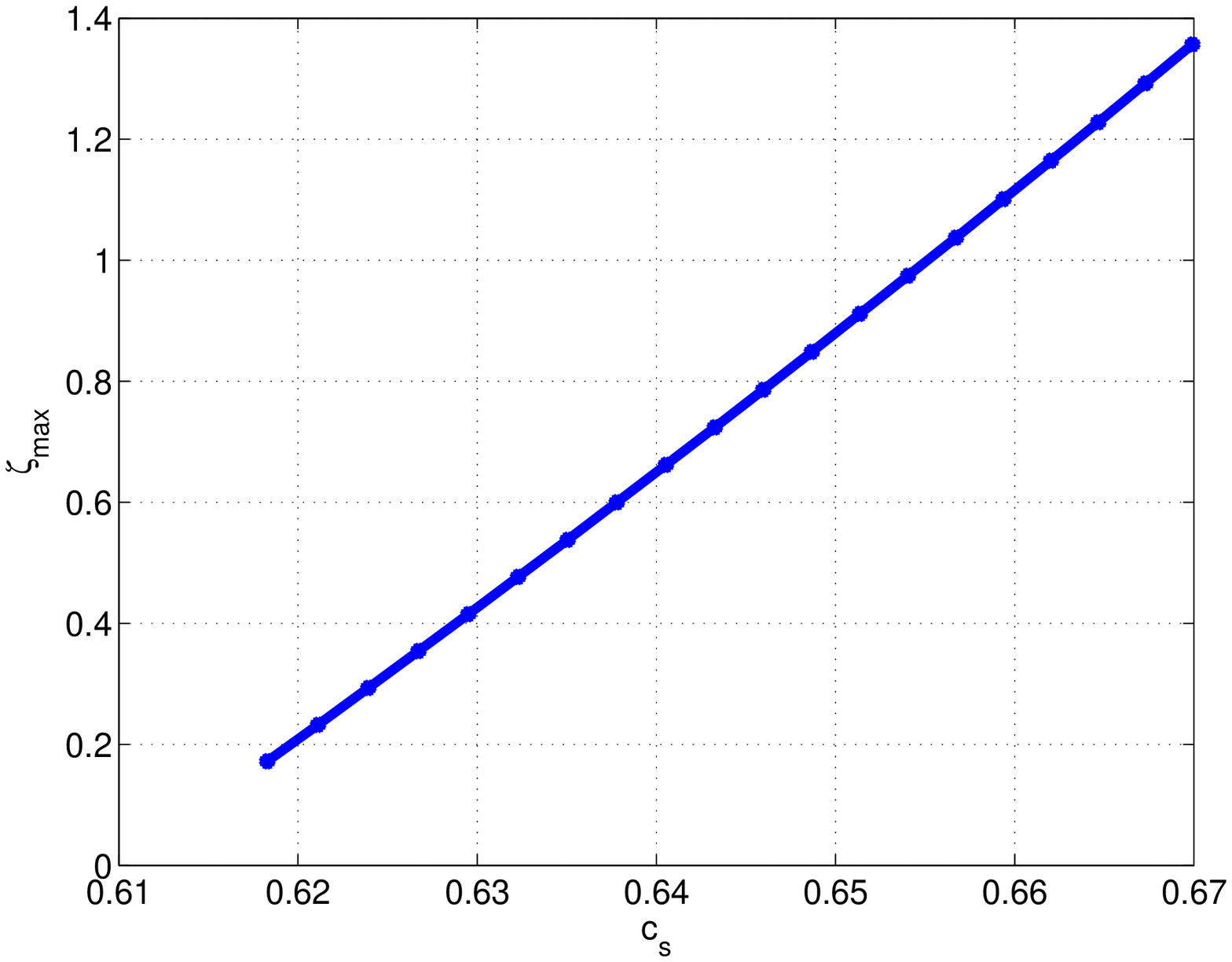}}
\subfigure[]
{\includegraphics[height=5.9cm,width=5.9cm]{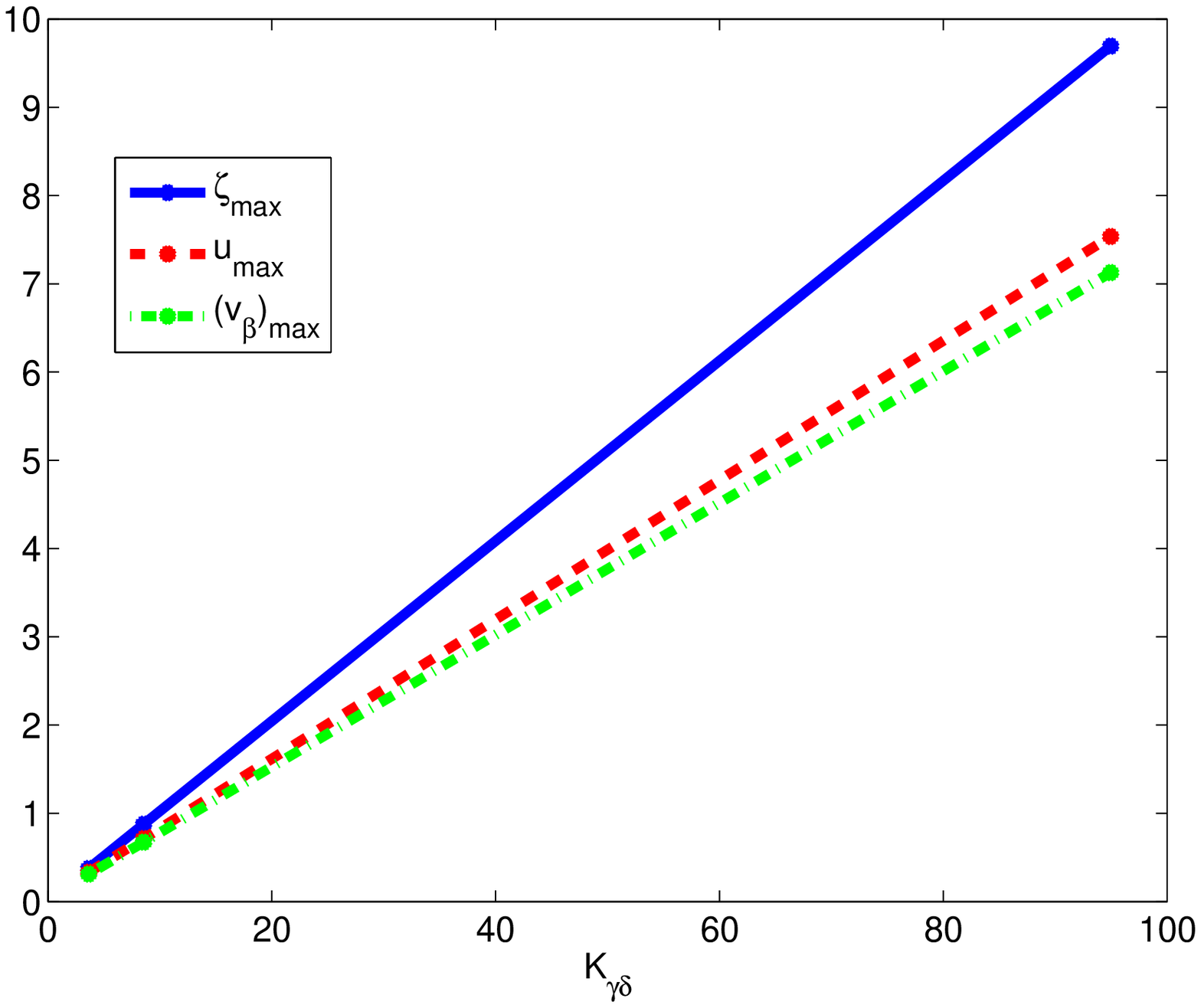}}
\caption{(a),(b) Amplitude vs. speed and (c) amplitude vs. $K_{\gamma,\delta}$ for the profiles of Figure \ref{f2}(a).}
\label{f3}
\end{figure}

This speed-amplitude relation is confirmed in Figure \ref{f3}, which corresponds to the waves of elevation of Figure \ref{f2}(a). In Figure \ref{f3}(a) we can observe the amplitude as increasing function of the speed in the deviation of the interface and both velocity variables. Figure \ref{f3}(b) shows the speed-amplitude relation only for the interfacial deviation $\zeta$ and for a wider range of speeds. The goal is trying to specify a relation $\zeta_{max}\approx Ac_{s}^{B}+C$ by fitting the experimental data for the parameters $A, B$ and $C$. The resolution of the fitting problem gives $A=18, B=2.75$ and $C=-4.626$ with, as a measure of the goodness of fit, a sum of squares due to error (SSE) of about $4.146\times 10^{-9}$, a R-squared value of $1$ and a root mean squared error (RMSE) of approximately $2.434\times 10^{-5}$.

\begin{figure}[htbp]
\centering
\subfigure[]
{\includegraphics[height=8.35cm,width=8.35cm]{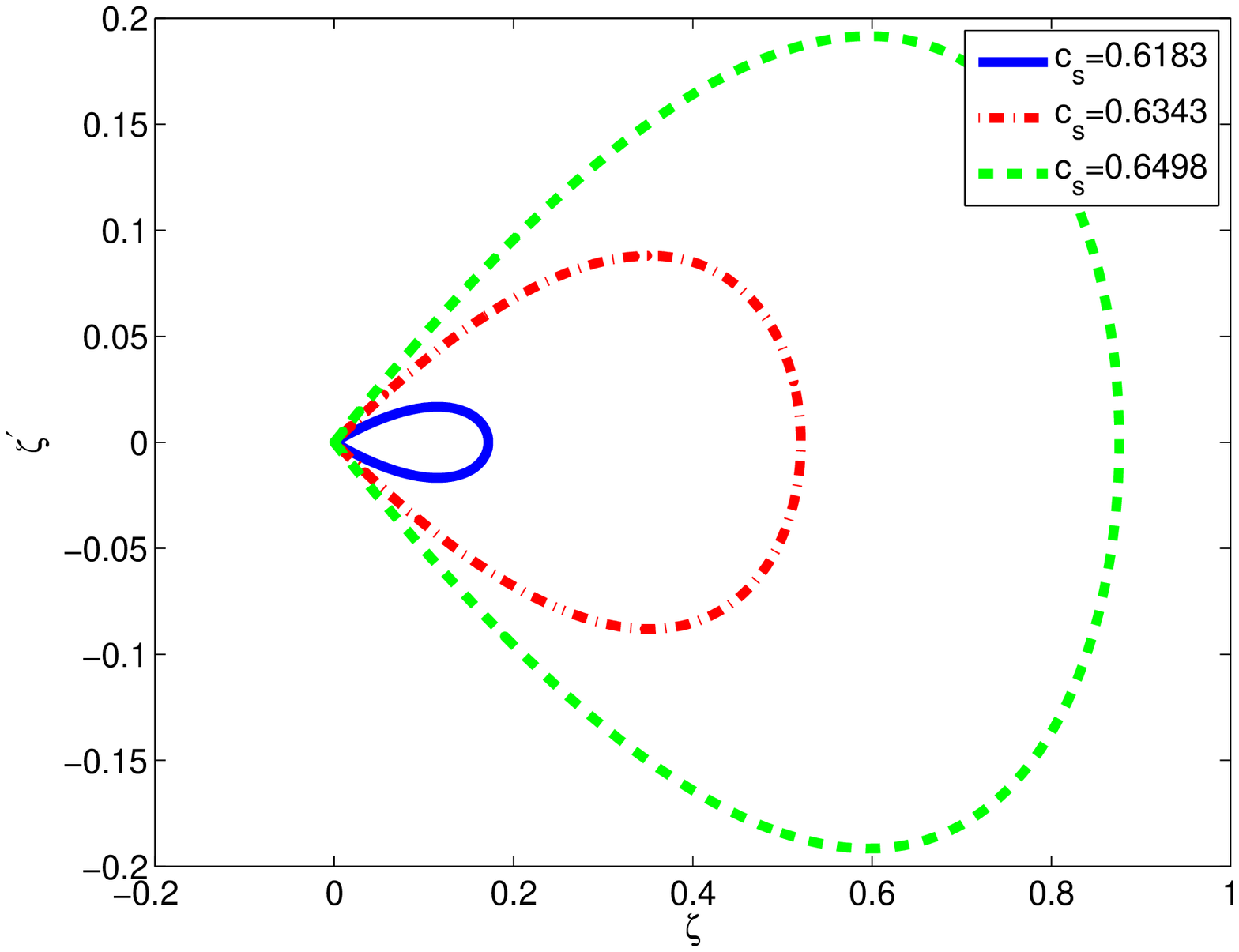}}
\subfigure[]
{\includegraphics[height=8.35cm,width=8.35cm]{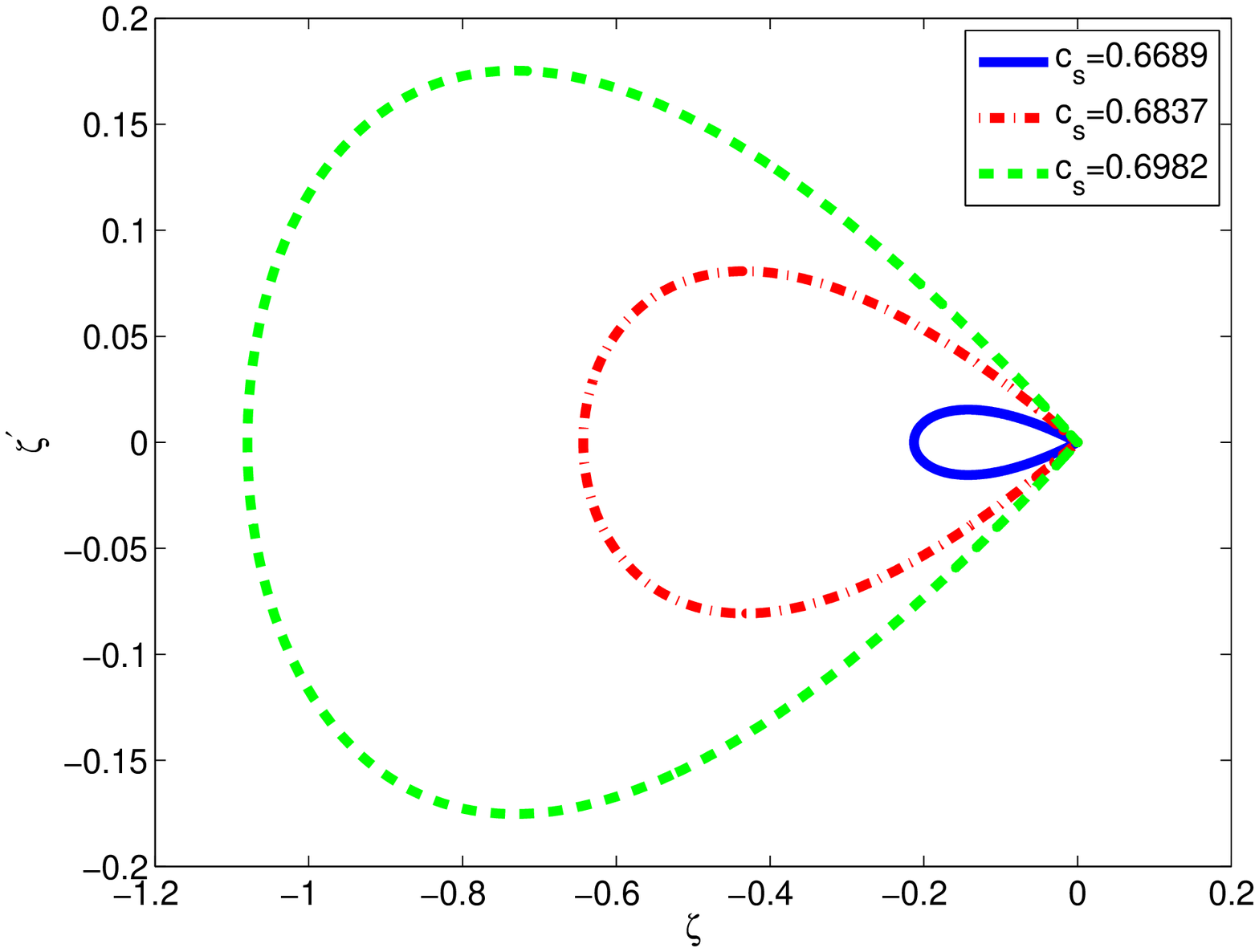}}
\caption{Phase portraits of the approximate profiles shown in Figure \ref{f2}.}
\label{f4}
\end{figure}

On the other hand, 
in Figure \ref{f3}(c), for $\gamma=0.5$, various $\delta$ and $c_{s}=c_{\gamma,\delta}+0.05$, the amplitude is displayed as function of the parameter $K_{\gamma,\delta}$ which determines the type (of elevation or of depression) of the profile. The amplitude also increases with $K_{\gamma,\delta}$ and Figure \ref{f3}(c) shows the dependence of the amplitude of the waves on the depth and density ratios of the two-fluid system (through the parameter $K_{\gamma,\delta}$).

\begin{figure}[htbp]
\centering
\subfigure[]
{\includegraphics[height=8.35cm,width=8.35cm]{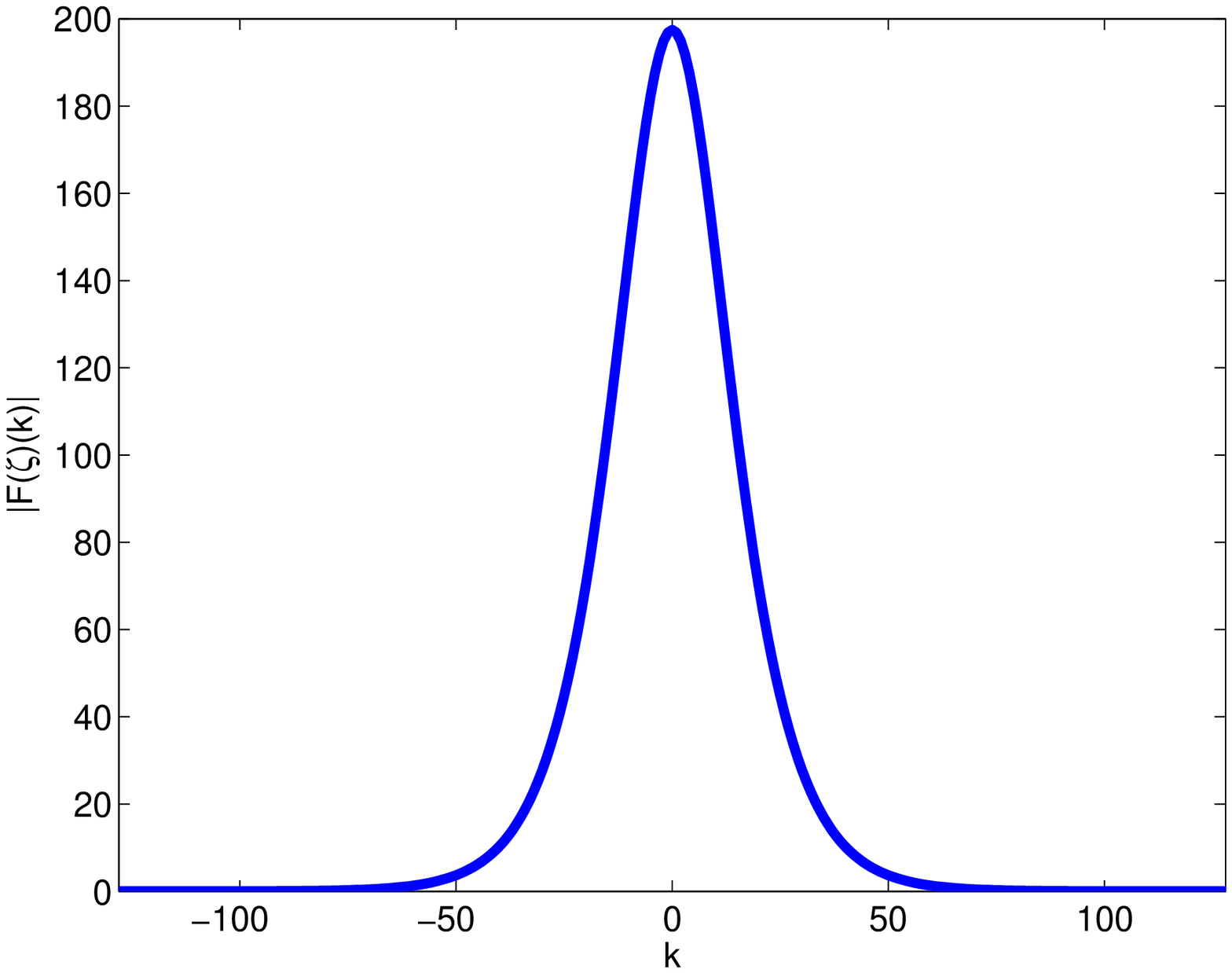}}
\subfigure[]
{\includegraphics[height=8.35cm,width=8.35cm]{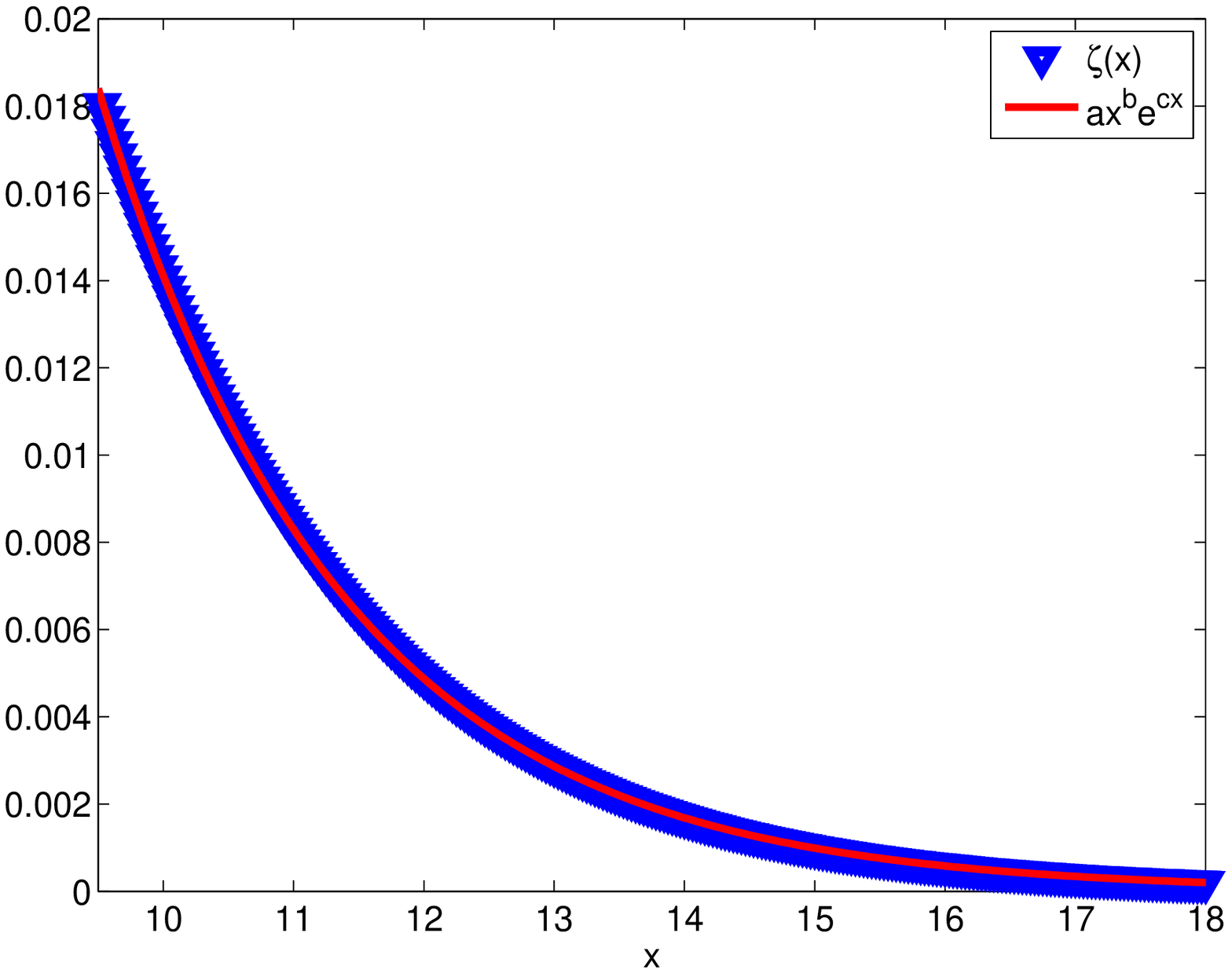}}
\caption{(a) Modulus of the Fourier transform $|\widehat{\zeta}(k)|$ of a profile vs. $k$; (b) Profile $\zeta(x)$ and fitting curve $ax^{b}e^{cx}$ vs. $x$.}
\label{f5}
\end{figure}

Figure \ref{f4} shows the phase portraits of the approximate profiles of Figure \ref{f2}, computed by using pseudospectral differentiation, \cite{Canutohqz}. We think that the smooth way how the orbits approach the origin suggests an exponential decay of the waves. This is also supported by the smooth form of the modulus of the Fourier transform of one of the profiles, displayed in Figure \ref{f5}(a). If part of the wave is fitted to a curve of the form $ax^{b}e^{cx}$ (Figure \ref{f5}(b)), then the solution of the corresponding least squares problem confirms this behaviour (see Table \ref{tav1}). Actually, similar arguments to those of \cite{PegoW1997} for the case of surface water waves (by applying the Stable Manifold Theorem) might be used to show that the profiles decay exponentially to zero. 

\begin{figure}[htbp]
\centering
\subfigure[]
{\includegraphics[height=8.35cm,width=8.35cm]{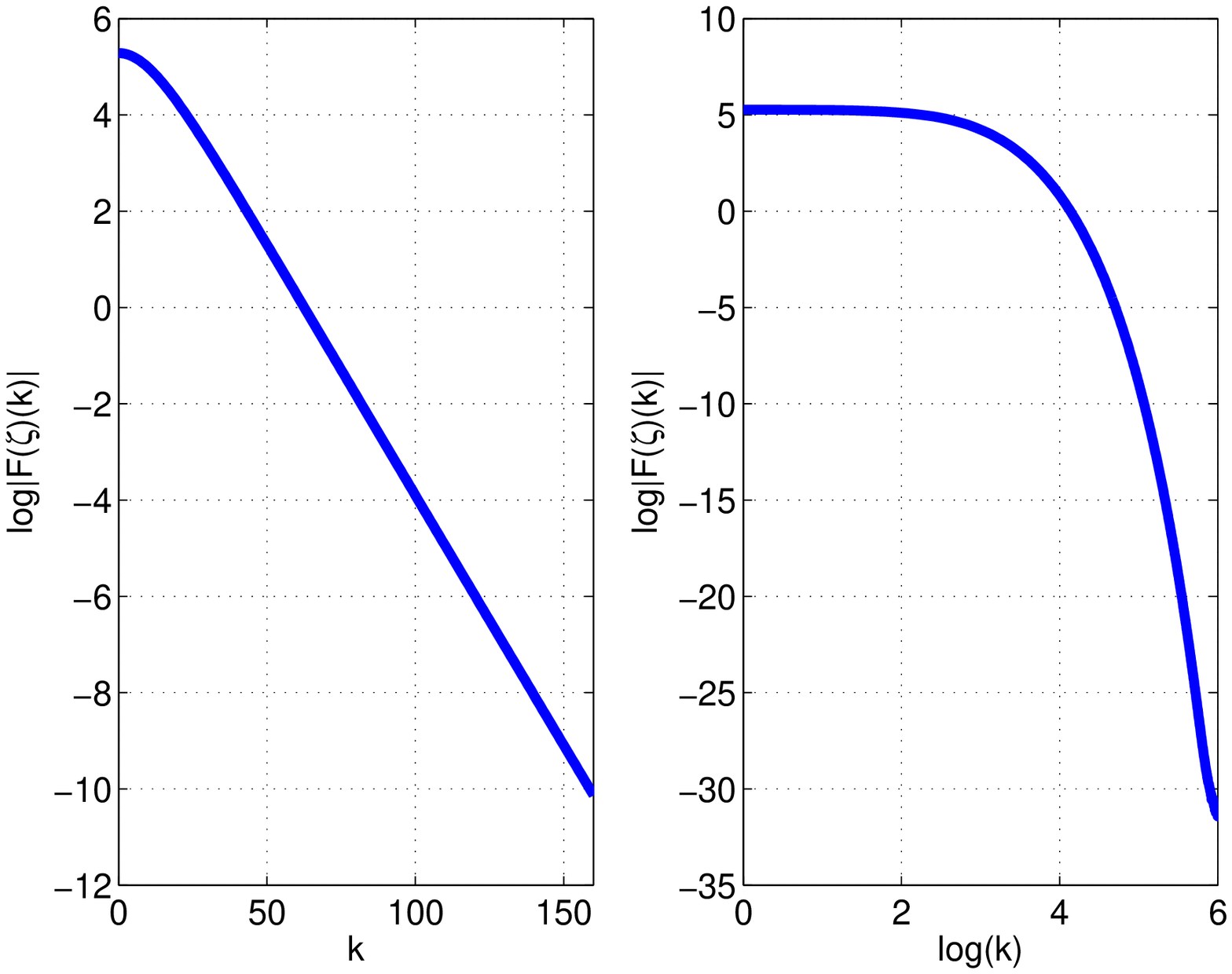}}
\subfigure[]
{\includegraphics[height=8.35cm,width=8.35cm]{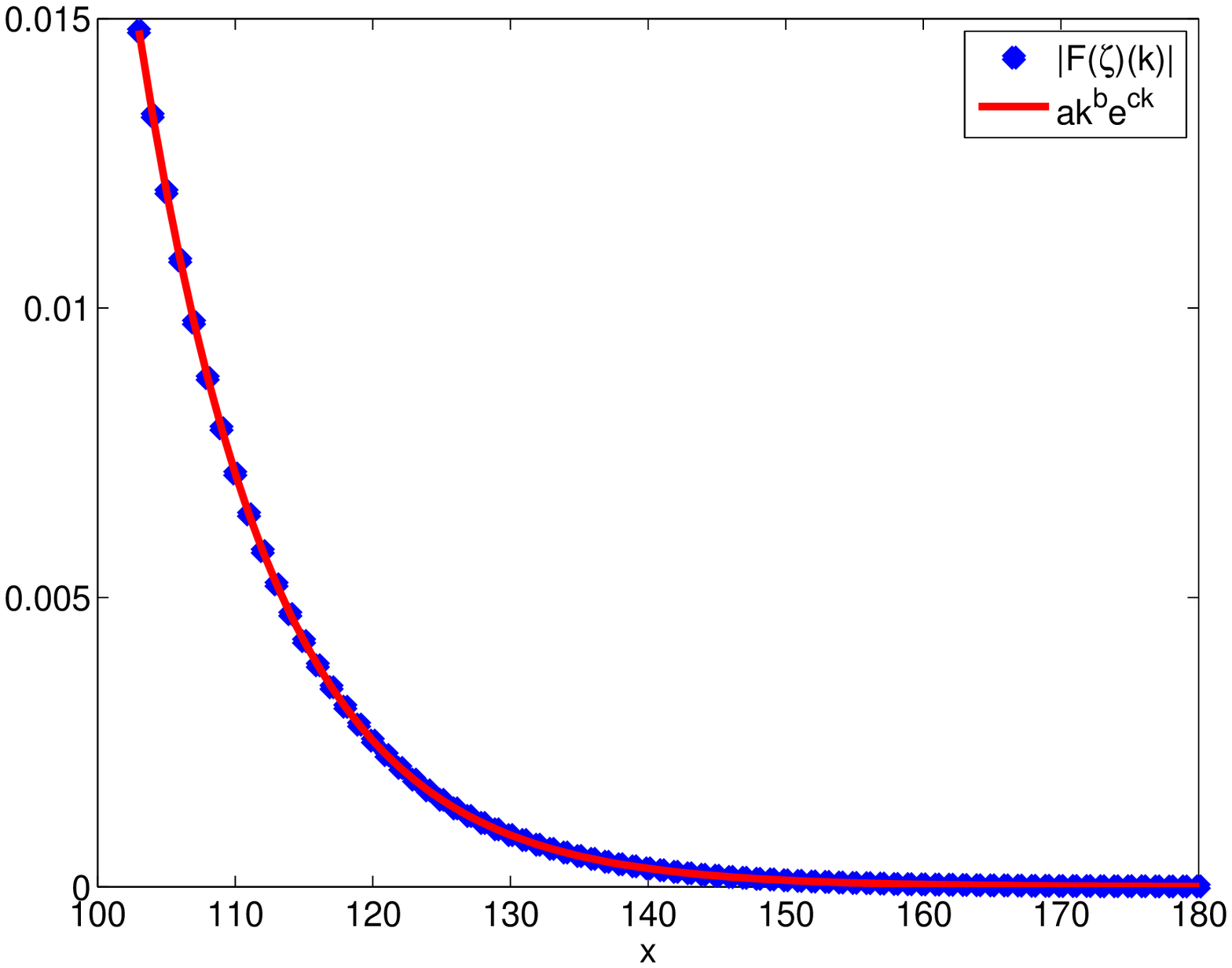}}
\caption{(a) Modulus of the Fourier transform vs. $k$ in linear-log and log-log scales, \cite{SulemSF1983}; (b) Modulus of the Fourier transform of a profile and fitting curve $ak^{b}e^{ck}$ vs. $k$..}
\label{f6}
\end{figure}

The decay of the Fourier transform can give more information on the regularity of the profiles. This is illustrated in Figure \ref{f6}. In linear-log and log-log scales, Figure \ref{f6}(a) shows the modulus of the Fourier transform versus the wavenumber, while in Figure \ref{f6}(b), part of this modulus is fitted to a curve of the form $ak^{b}e^{ck}$, \cite{SulemSF1983}, and according to the goodness of fit (Table \ref{tav1}), the curve shows a reliable exponential decay, which suggests a smooth character of the solitary wave profile.

\begin{table}[htbp]
\begin{center}
  \begin{tabular}{|c|c|}
    \hline\hline
    {\textit{fitting curve}}  &{\textit{g.o.f.}}
    \\
    \hline
$f(x)=ax^{b}e^{cx}$&$SSE=1.967\times 10^{-8}$\\  
$a=2.8176, b=0.0114$&$R-squared=1$\\  
$c=-0.5323$&$RMSE=4.233\times 10^{-6}$\\  
  \hline\hline
 $f(k)=ak^{b}e^{ck}$&$SSE=3.548\times 10^{-11}$\\  
$a=1261.4, b=-0.1728$&$R-squared=1$\\  
$c=-0.1023$&$RMSE=4.896\times 10^{-7}$\\ 
\hline\hline  
  \end{tabular}
  \caption{Curves and goodness of fit for data from Figures \ref{f5}(b) (up) and \ref{f6}(b) (down).}
  \label{tav1}
  \end{center}
\end{table}
\section{Conclusions}
\label{sec:sec6}
Considered in this paper is a nonlocal system for internal waves. The model was derived in \cite{BLS2008} as a consistent approximation to the Euler equations for the wave propagation in a two-layer fluid system, with rigid lid assumptions and under a Boussinesq physical regime for both fluids. In this paper, having in mind the treatment of the case of surface wave propagation in the literature, \cite{BChS2002,PegoW1997}, several mathematical aspects of the model are studied in more detail. After an alternative derivation of the differential system (which does not make use of the nonlocal operators considered in the original presentation in \cite{BLS2008}), a result of local existence and uniqueness of solution in suitable Sobolev spaces is proved. The study is then focused on the solitary wave solutions. They are shown to exist for a range of the speeds which includes bidirectional propagation and depends on the depth and density ratios of the two-fluid problem. These ratios also determine the character of elevation or of depression of the wave. Finally, a numerical tecnique to generate approximate solitary wave profiles is introduced and applied to suggest some features: the waves look to be smooth, decay exponentially to zero at infinity and their amplitude looks to behave as a power function of the speed.

The analysis, theoretical and numerical, of the solitary waves performed in the last part of the paper must be taken into account as starting point for a future research, mainly focused on the dynamics of the waves. In this sense, a computational work concerning experiments of stability under small and large perturbations, head-on and overtaking collisions, resolution of initial data into trains of solitary waves, etc, is essential for this purpose. This includes a detailed analysis of the codes designed and used to this end. From the theoretical point of view, the results, presented in the literature (see e.~g. \cite{PegoW1997}) about stability are again a necessary tool to go more deeply into this topic.

\section*{Acknowledgements}
This research is supported by Spanish Ministerio de Econom\'{\i}a y Competitividad under grant MTM2014-54710-P with the participation of FEDER. The author would like to thank Professors V. Dougalis and D. Mitsotakis for fruitful discussions and so important suggestions.



\end{document}